\documentclass[twoside]{article}
\usepackage{landau}
\usepackage[dvips]{epsfig}
\theoremstyle{definition}

\theoremstyle{remark}

\numberwithin{equation}{section}

\def\shorttitle{Computing material fronts with a Lagrange-Projection approach}
\def\shortauthor{C. Chalons, F. Coquel}

\begin{document}

\title
{\Large \bf \boldmath Computing material fronts with a Lagrange-Projection approach\footnote{Research 
of the first author was partially supported by the French Agence Nationale de la Recherche, contract ANR-08-JCJC-0132-01. 
Research of the authors was partially supported by the Minist\`ere de la Recherche under grant ERTint 1058 entitled 
{\it Simulation Avanc\'ee du Transport des Hydrocarbures.}}}%

\author{\large Christophe Chalons
    \\ \normalsize\emph{CEA-Saclay, DEN/DANS/DM2S/SFME/LETR}
    \\ \normalsize\emph{F-91191 Gif-sur-Yvette, France.}
    \\ \normalsize\emph{E-mail: chalons@math.jussieu.fr}
\\[2mm]
    \large Fr\'ed\'eric Coquel
    \\ \normalsize\emph{CNRS \& Centre de Math\'ematiques Appliqu\'ees, U.M.R.~7641}
    \\ \normalsize\emph{Ecole Polytechnique, Route de Saclay, 91128 Palaiseau Cedex, France.}
    \\ \normalsize\emph{E-mail: frederic.coquel@cmap.polytechnique.fr}
}

\date{\vspace{-12mm}}

\maketitle

\thispagestyle{first}

\begin{abstract}
This paper reports investigations on the computation of material fronts in multi-fluid 
models using a Lagrange-Projection approach. Various forms of the Projection step are considered.
Particular attention is paid to minimization of conservation errors.
\end{abstract}

\section{Introduction}
It is well-known that standard conservative discretizations
of gas dynamics equations in Eulerian coordinates generally develop non physical pressure oscillations 
near contact discontinuities, and 
more generally near material fronts in multi-component flows. Several cures based on a local non conservative modification 
have been proposed. Let us quote for instance the hybrid algorithm 
derived by Karni in \cite{karni2} and the Two-Flux method proposed by Abgrall and Karni
in \cite{ak1} for multi-fluid flows. See also \cite{abgrall2}, \cite{bhr1} and the references therein. \\
We investigate here a Lagrange-Projection type
method to get rid of pressure oscillations. The basic motivation lies in the fact that oscillations do not 
exist in Lagrangian computations. It is then possible to clearly determine which operation in the projection 
step sparks off pressure oscillations. As in \cite{karni2}, \cite{ak1}, a non conservative correction is proposed. 
It is based on a local pressure averaging and a random sampling strategy on the mass fraction in order 
to strictly preserve isolated material fronts and get a statistical conservation property. Numerical results 
are proposed and compared with the 
Two-Flux method \cite{ak1}.

\section{The model under consideration}

We consider a nonlinear partial differential equations model governing the flow of two species 
$\Sigma_1$ and $\Sigma_2$ separated by a material interface. For instance, we focus on two perfect gases 
and we set
$p_{i} (\rho_{i},e_i) = (\gamma_{i} - 1) \rho_{i} e_i$, $\gamma_i = {C_{p,i}}/{C_{v,i}}$
and $T_{i} (\rho_{i},e_i) = {e_i}/C_{v,i}$
where $p_{i}$, $\rho_{i}$, $e_{i}$, $\gamma_{i} > 1$, $T_{i}$, $C_{p,i}>0$, $C_{v,i}>0$ respectively denote the pressure, the density, the internal energy, the adiabatic coefficient, the temperature  and the specific heats of 
$\Sigma_i$, $i=1,2$. The mixture density is given by $\rho = \rho_1+\rho_2$ and we adopt a Dalton's law for the mixture 
pressure $p= p_{1} (\rho_{1},e_1)+p_{2} (\rho_{2},e_2)$. 
We assume in addition that the two species evolve according to the same velocity $u$ and are at thermal equilibrium, that is 
$T=T_1(\rho_1,e_1)=T_2(\rho_2,e_2)$.
The mixture internal and total energies $e$ and $E$ are defined by 
$\rho e = \rho_1 e_1 + \rho_2 e_2$ and $\rho E = \frac{1}{2} \rho u^2 + \rho e$. Then, 
introducing the mass fraction $Y = \rho_1 / \rho$, straightforward manipulations yield 
$\rho e = \rho C_v T$ with $C_v = C_v(Y) = Y C_{v,1} + 
(1 - Y) C_{v,2}$ and 
$$
p = (\gamma - 1) \rho e \quad \mbox{with} \quad  
\gamma = \gamma(Y) = \frac{Y C_{p,1} + (1-Y) C_{p,2}}{Y C_{v,1} + (1-Y) C_{v,2}} > 1.
$$
In one-space dimension, the model under consideration writes 
\begin{equation} \label{systeme}
\left\{
\begin{array}{l}
\partial_t \rho + \partial_x (\rho u) = 0,\\
\partial_t (\rho u) + \partial_x (\rho u^2 + p) = 0,\\
\partial_t (\rho E) + \partial_x (\rho E + p)u = 0, \\
\partial_t \rho Y + \partial_x (\rho Y u) = 0,
\end{array}
\right.
\end{equation}
and for the sake of conciseness, we set
\begin{equation} \label{systemec}
\partial_t {\bf u}(x,t) + \partial_x {\bf f}({\bf u}(x,t)) = 0. \\
\end{equation}
The flux function ${\bf f}$ finds a natural definition with respect to 
the conservative unknowns ${\bf u} = (\rho, \rho u, \rho E, \rho Y)$.
Let us mention that (\ref{systemec}) is hyperbolic with eigenvalues 
$\lambda_0({\bf u}) = u$ and $
\lambda_{\pm}({\bf u}) = u \pm c({\bf u})$, 
$c({\bf u}) = \sqrt{{\gamma p}/{\rho}}$, provided that $\rho >0$, $0 \leq Y \leq 1$ and $p>0$. 
The characteristic field 
associated with $\lambda_0$ is linearly degenerate, leading to 
{\it contact discontinuities} or {\it material fronts}. The two extreme fields are genuinely nonlinear.

\section{Numerical schemes}
 
This section is devoted to the discretization of (\ref{systemec}). As already stated, 
a specific attention must be paid to the contact discontinuities to avoid 
pressure oscillations. With this in mind, we first revisit the "Two-Flux Method" proposed by Abgrall and Karni \cite{ak1} and then 
present a new numerical procedure based on a Lagrangian approach and a random sampling strategy. Comparisons 
will be proposed in section \ref{sec:numexp}. \\
Let us introduce a time step $\Delta t>0$ and a space step 
$\Delta x>0$ that we assume to be constant for simplicity. We set $\lambda = {\Delta t}/{\Delta x}$ and define
the mesh interfaces $x_{j+1/2}=j \Delta x$ for $j \in \mathbb{Z}$, and the 
intermediate times $t^n=n \Delta t$ for $n \in \mathbb{N}$. In the sequel,
${\bf u}^n_j$ denotes the approximate value of ${\bf u}$ at time 
$t^n$ and on the cell $\mathcal{C}_j = [x_{j-1/2}, x_{j+1/2}[$.
For $n=0$ and $j \in \mathbb{Z}$, we set 
$
{\bf u}^0_j = \frac{1}{\Delta x} \int_{x_{j-1/2}}^{x_{j+1/2}} {\bf u}_0(x) dx
$ 
where ${\bf u}_0(x)$ is the initial condition.

\subsection{The Two-Flux method revisited} \label{ttfm}

Aim of this section is to review the Two-Flux method proposed by Abgrall and Karni \cite{ak1}. 
Let us first recall that pressure oscillations 
do not systematically appear in 
{\it single-fluid} computations. Abgrall and Karni \cite{ak1} then propose to replace any conservative 
{\it multi-fluid} strategy by a non conservative approach based on the definition 
of two {\it single-fluid} numerical fluxes at each interface. 
We first 
recall the algorithm in details and then suggest a slight modification in order 
to lessen the conservation errors. This strategy will be used as a reference to assess the validity 
of the Lagrangian strategies proposed in the next subsection. \\
\ \\
{\bf The original algorithm.} Let us consider a two-point numerical flux function 
${\bf g}$ consistent with ${\bf f}$. 
The Two-Flux method proposes to update the 
sequence $({\bf u}^n_j)_{j \in \mathbb{Z}}$ into two steps, under an usual $1/2$ CFL 
restriction. \\
\ \\
{\it First step : evolution of $\rho$, $\rho u$ and $p$} ($t^n \to t^{n+1-}$) \\
Let us first define ${\bf v}= (\rho, \rho u, p, Y)$ and the one-to-one mapping
${\bf u} \to {\bf v} = {\bf v}({\bf {u}})$ thanks to the thermodynamics closures.
Two interfacial numerical fluxes ${\bf g}_{j+1/2,L}$ and
${\bf g}_{j+1/2,R}$ are then defined by
\begin{equation}
{\bf g}_{j+1/2,L} = {\bf g}({\bf u}^{n}_{j},{\bf \overline{u}}^{n}_{j+1,L}), \quad
{\bf g}_{j+1/2,R} = {\bf g}({\bf \overline{u}}^{n}_{j,R},{\bf u}^{n}_{j+1}),
\end{equation}
where ${\bf \overline{u}}^{n}_{j+1,L} = {\bf u}({\bf \overline{v}}^{n}_{j+1,L})$ with
${\bf \overline{v}}^{n}_{j+1,L} = (\rho^{n}_{j+1},(\rho u)^{n}_{j+1},
p^{n}_{j+1}, Y^n_j),$
and
${\bf \overline{u}}^{n}_{j,R} = {\bf u}({\bf \overline{v}}^{n}_{j,R})$ with
${\bf \overline{v}}^{n}_{j,R} = (\rho^{n}_{j},(\rho u)^{n}_{j},
p^{n}_{j}, Y^n_{j+1})$. In some sense, the mass fraction $Y$ is then 
assumed to be the same on both side 
of each interface since 
$Y=Y_{j}^n$ is used for the computation of ${\bf g}_{j+1/2,L}$ and
$Y=Y_{j+1}^n$ for ${\bf g}_{j+1/2,R}$. At last, we use
${\bf g}_{j+1/2,L}$, respectively ${\bf g}_{j+1/2,R}$, to update the conservative 
unknowns $\rho$, $\rho u$ and $\rho E$ on the cell $j$, resp. $(j+1)$. 
With clear notations, we get for all
$j \in \mathbb{Z}$
\begin{equation} \label{mgfs1a}
\begin{array}{l} 
{\bf \rho}^{n+1-}_{j} = 
{\bf \rho}^n_{j} - 
  \lambda 
  ({\bf g}_{j+1/2,L}^{\rho} - 
  {\bf g}_{j-1/2,R}^{\rho}),\\
{(\rho u)}^{n+1-}_{j} = 
{(\rho u)}^n_{j} - 
  \lambda 
  ({\bf g}_{j+1/2,L}^{\rho u} - 
  {\bf g}_{j-1/2,R}^{\rho u}),\\
{(\rho E)}^{n+1-}_{j} = 
{(\rho E)}^n_{j} - 
  \lambda 
  ({\bf g}_{j+1/2,L}^{\rho E} - 
  {\bf g}_{j-1/2,R}^{\rho E}).
\end{array}
\end{equation}
Let us note that $\rho_1=\rho Y$ is not concerned with 
(\ref{mgfs1a}). We simply set
${(\rho_1)}^{n+1-}_{j} = {\rho}^{n+1-}_{j} \times
{Y}^n_{j}$
($Y$ is naturally kept constant in this step) and then
$p^{n+1-}_{j} = p({\bf u}^{n+1-}_{j})$ with 
${\bf u}^{n+1-}_{j} = (\rho,\rho u,\rho E,\rho_1)^{n+1-}_{j}.$
\ \\
\ \\
{\it Second step : evolution of $\rho_1$} ($t^{n+1-} \to t^{n+1}$) \\
In this step, $\rho_1$ is evolved in a conservative way using the numerical flux function 
$g$ while the values of 
$\rho$, $\rho u$ and $p$ are kept unchanged. Again with clear notations, 
the vector
${\bf v}^{n+1}_{j} = (\rho^{n+1}_{j},(\rho u)^{n+1}_{j},p^{n+1}_{j},Y^{n+1}_{j})$ 
is then defined 
${\bf \rho}^{n+1}_{j} = 
{\bf \rho}^{n+1-}_{j}$, 
${(\rho u)}^{n+1}_{j} = 
{(\rho u)}^{n+1-}_{j}$, 
${p}^{n+1}_{j} =
p^{n+1-}_{j}$
and $Y^{n+1}_{j} = (\rho_1)^{n+1}_{j} / \rho^{n+1}_{j}$, where 
\begin{equation}
{(\rho_1)}^{n+1}_{j} = 
{(\rho_1)}^{n}_{j} - 
  \lambda 
  ({\bf g}_{j+1/2}^{\rho Y} - 
  {\bf g}_{j-1/2}^{\rho Y}),
\end{equation}  
with ${\bf g}_{j+1/2}={\bf g}({\bf u}^n_{j},{\bf u}^n_{j+1})$. At last, 
we set ${\bf u}^{n+1}_{j} = {\bf u}({\bf v}^{n+1}_{j})$. \\
\ \\
It is worth noticing that the Two-Flux method is not conservative on $\rho$ and $\rho u$ since the fluxes 
${\bf g}_{j+1/2,L}$ and ${\bf g}_{j+1/2,R}$ defined at 
each interface are different as soon as 
${Y}^{n}_{j} \neq {Y}^{n}_{j+1}$. It is not conservative on $\rho E$ either, but by construction the Two-Flux method is 
conservative on $\rho_1$. \\
\ \\
{\bf The associated quasi-conservative algorithm.} 
It is actually clear from \cite{karni2}, \cite{ak1} and the references therein that 
in standard conservative discretizations of 
(\ref{systemec}), only 
the update formula of the total energy $\rho E$ is responsible for the pressure oscillations. 
We are then tempted to propose a quasi-conservative variant of the Two-Flux method such that 
only the total energy is treated in a non conservative form.
For all $j \in \mathbb{Z}$, we simply replace 
(\ref{mgfs1a}) by
\begin{equation} \label{mgfs1abis}
\begin{array}{l} 
{\bf \rho}^{n+1-}_{j} = 
{\bf \rho}^n_{j} - 
  \lambda 
  ({\bf g}_{j+1/2}^{\rho} - 
  {\bf g}_{j-1/2}^{\rho}),\\
{(\rho u)}^{n+1-}_{j} = 
{(\rho u)}^n_{j} - 
  \lambda 
  ({\bf g}_{j+1/2}^{\rho u} - 
  {\bf g}_{j-1/2}^{\rho u}),\\
{(\rho E)}^{n+1-}_{j} = 
{(\rho E)}^n_{j} - 
  \lambda 
  ({\bf g}_{j+1/2,L}^{\rho E} - 
  {\bf g}_{j-1/2,R}^{\rho E}),
\end{array}
\end{equation}
where again ${\bf g}_{j+1/2}={\bf g}({\bf u}^n_{j},{\bf u}^n_{j+1})$. 
The second step is 
unchanged. 

\subsection{The Lagrange-Projection approach with random sampling}

In this section, we propose a Lagrangian approach for approximating the solutions of (\ref{systemec}). The general idea 
is to first solve this system in Lagrangian coordinates, and then to come back to an Eulerian description of the flow 
with a projection step. Under its classical conservative form, the Lagrange-Projection method 
generates spurious oscillations near the material fronts. In order to remove these oscillations, we 
propose to adapt the projection step (only). We begin with a description of the Lagrangian step and then recall, for the sake of clarity, the usual conservative projection step (see for instance \cite{gr2}). Again, an usual $1/2$ CFL restriction is used.\\
\ \\
{\it{The Lagrangian step ($t^n \rightarrow t^{n+1-}$)}} \\
In this step, (\ref{systemec}) in written in Lagrangian coordinates and solved 
by an acoustic scheme (see for instance \cite{despres1}), which gives
\begin{equation} \label{acoussch}
\begin{array}{l}
{{
{\tau}^{n+1-}_{j} = 
{\tau}^n_{j} - 
\lambda \tau^n_{j} 
(u^{n}_{j+1/2} - u^{n}_{j-1/2}), 
}}\\
{{
{u}^{n+1-}_{j} = 
{u}^n_{j} - 
\lambda \tau^n_{j} 
(p^{n}_{j+1/2} - p^{n}_{j-1/2}), 
}}\\
{{
{E}^{n+1-}_{j} = 
{E}^n_{j} - 
\lambda \tau^n_{j}
\big((pu)^{n}_{j+1/2} - (pu)^{n}_{j-1/2}\big),  
}}
\\
{{
{Y}^{n+1-}_{j} = 
{Y}^n_{j},  
}}
\end{array}
\end{equation}
where the velocity and the pressure at interfaces are defined by
\begin{equation} \label{defupintlp}
\left\{
\begin{array}{l}
{\displaystyle{
u^{n}_{j+1/2} = \frac{1}{2}(u^n_{j}+u^n_{j+1}) + 
\frac{1}{2(\rho c)^n_{j+1/2}}(p^n_{j}-p^n_{j+1}), 
}}\\
{\displaystyle{
p^{n}_{j+1/2} = \frac{1}{2}(p^n_{j}+p^n_{j+1}) + 
\frac{(\rho c)^n_{j+1/2}}{2}(u^n_{j}-u^n_{j+1}). 
}}
\end{array}
\right.
\end{equation}
The proposed local approximation $(\rho c)^n_{j+1/2}$ of the Lagrangian sound speed is
$(\rho c)^n_{j+1/2} = \max((\rho c)^n_{j}, (\rho c)^n_{j+1})$ 
but other definitions may be found for instance in \cite{lagoutiere1}.
In this step, the grid points $x_{j+1/2}$ move at velocity $u^n_{j+1/2}$ so that 
$\rho^{n+1-}_{j} =1/\tau^{n+1-}_{j}$,
$(\rho u)^{n+1-}_{j} = \rho^{n+1-}_{j} \times u^{n+1-}_{j}$,
$(\rho E)^{n+1-}_{j} = \rho^{n+1-}_{j} \times E^{n+1-}_{j}$
and
$(\rho_1)^{n+1-}_{j} = \rho^{n+1-}_{j} \times Y^{n+1-}_{j}$
define approximate values of ${\bf u}$
on a Lagrangian grid with mesh interfaces 
$x_{j+1/2}^* = x_{j+1/2} + u^n_{j+1/2} \Delta t$. \\
\ \\
{\it{The usual projection step ($t^{n+1-} \rightarrow t^{n+1}$)}} \\
Aim of this step is to project the solution obtained at the end of the first step on the Eulerian grid 
defined by the mesh interfaces $x_{j+1/2}$. {\it Usually}, the choice is made to project the conservative vector 
${\bf u}$ in order to obtain a conservative Lagrange-Projection scheme (see again \cite{gr2}). 
More precisely, such a choice writes  
\begin{equation} \label{moyL2phi1}
\varphi^{n+1}_{j}  = \frac{1}{\Delta x} 
\int_{x_{j-1/2}}^{x_{j+1/2}} \varphi^{n+1-}(x)dx \quad \mbox{with} \quad
\varphi=\rho, \rho u, \rho E, \rho_1=\rho Y,
\end{equation}
or, with
$\Delta x_{j}^* = x_{j+1/2}^* - x_{j-1/2}^*$ and 
$\varepsilon(j,n)=
\left\{
\begin{array}{rcl}
-1/2 & \mbox{if} & u^n_{j+1/2} >0, \\
1/2 & \mbox{if} & u^n_{j+1/2} <0, 
\end{array}
\right.
$, 
$$
\varphi^{n+1}_{j} = \frac{1}{\Delta x} 
\{
\Delta x_{j}^* \varphi^{n+1-}_{j} - 
\Delta t (u^n_{j+1/2} \varphi^{n+1-}_{j+1/2+\varepsilon(j,n)} -
u^n_{j-1/2} \varphi^{n+1-}_{j-1/2+\varepsilon(j-1,n)}) 
\}.
$$
{\bf What is wrong with this scheme ?} 
The Lagrange-Projection approach allows to precisely reveal the operation that makes 
the material fronts necessarily damaged by pressure oscillations.
Let us indeed consider an isolated material front with uniform velocity and pressure 
profiles~: $u(x,0) = u^0 > 0$ and $p(x,0) = p^0$, while $Y(x,0) = 1$ if $x<0$ and $Y(x,0) = 0$ if $x>0$. 
The density is also set to be uniform for simplicity :  
$\rho(x,0) = \rho^0$. We first observe that this profile is clearly preserved in the Lagrangian step since 
by (\ref{defupintlp})
we have $u^0_{j+1/2} = u^0$ and $p^0_{j+1/2} = p^0$ by (\ref{defupintlp}). Then, 
the projection procedure (\ref{moyL2phi1}) gives $\rho^{1}_j = \rho^0$ and 
$(\rho u)^{1}_j = \rho^0 u^0$
so that
$u^1_j = (\rho u)^{1}_j/\rho^{1}_j = u^0$.
After the first time iteration, the velocity profile is then still free of spurious oscillations. 
At last, using the property that this velocity is constant and positive,
and focusing for instance on the cell of index $j=1$, (\ref{moyL2phi1}) gives for $\varphi = \rho E, \rho_1$  
$$
(\rho E)^{1}_{1} = (\rho E)^0_1 -  
u^0 \lambda \big((\rho E)^{0}_{1} - (\rho E)^{0}_{0}\big) \quad \mbox{and} \quad Y_1^1 = u^0 \lambda. 
$$
The pressure is then given after easy calculations by
$$
p^1_1 = 
p^0 \times \big(\gamma(Y^1_1) - 1\big) \times \big(\frac{1}{\gamma_2-1}(1 - u^0 \lambda) + u^0 \lambda\frac{1}{\gamma_1-1}\big).
$$
At this stage, there is no reason for $p^1_1$ to equal $p^0$ 
From the very first time iteration, a pressure oscillation is 
then created. As an immediate consequence, the velocity profile will not remain uniform in the next time iteration
and the numerical solution is damaged for good. \\
\ \\
It is then clear that the way the pressure 
is updated in the projection step (only) is responsible for the spurious oscillations in an usual 
conservative Lagrange-Projection scheme. We propose to modify this step. 
As for the Two-Flux method, the idea is to give up the conservation property in order to maintain uniform the pressure $p$ across material fronts. Let us emphasize that the Lagrangian step is unchanged.\\
\ \\
{\it{The quasi-conservative $p$-projection step ($t^{n+1-} \rightarrow t^{n+1}$)}} \\
First of all, $\rho$, $\rho u$ and $\rho_1$ still evolve according to (\ref{moyL2phi1}) 
so that the algorithm remains conservative on these variables. 
We will keep on using (\ref{moyL2phi1}) for $\varphi=\rho E$ 
only for $j$ not in a subset $\mathbb{Z}_p^{\varepsilon}$ of $\mathbb{Z}$ defined below. On 
the contrary, the pressure $p$ (instead of $\rho E$) is 
averaged for $j \in \mathbb{Z}_p^{\varepsilon}$: 
\begin{equation} \label{newp}
p^{n+1}_{j}  = \frac{1}{\Delta x} 
\int_{x_{j-1/2}}^{x_{j+1/2}} p^{n+1-}(x)dx.
\end{equation}
For $j \in \mathbb{Z}_p^{\varepsilon}$, we then set 
${\bf v}^{n+1}_{j} = (\rho^{n+1}_{j},(\rho u)^{n+1}_{j},p^{n+1}_{j},Y^{n+1}_{j})$ with
$Y^{n+1}_{j} = (\rho_1)^{n+1}_{j} / \rho^{n+1}_{j}$, and 
${\bf u}^{n+1}_{j} = {\bf u}({\bf v}^{n+1}_{j})$. \\
\ \\
{\bf Definition of $\mathbb{Z}_p^{\varepsilon}$.} 
Up to our knowledge, the idea of averaging the pressure $p$ in a Lagrange-Projection strategy 
first appeared in \cite{bhr1}. 
This way to proceed 
is clearly sufficient to remove the pressure oscillations 
near the material fronts if ${{\mathbb{Z}_p^{\varepsilon}}}={{\mathbb{Z}}}$.
However, averaging the pressure $p$ 
instead of the total energy $\rho E$ for all $j \in \mathbb{Z}$ gives a non conservative scheme 
that is expected to 
provide discontinuous solutions violating the Rankine-Hugoniot conditions, see for instance 
Hou and LeFloch \cite{hl1} (note however that here, 
the Lagrangian system associated with (\ref{systeme}) is actually treated in
a {\it conservative} form, while the pressure averaging takes place in the projection step only). 
This was confirmed in practice when considering solutions involving 
shocks with large amplitude. 
Then, in order to lessen the conservation 
errors, we propose to localize the non conservative treatment around the contact discontinuities setting 
$
{{\mathbb{Z}_p^{\varepsilon}}} = \{{{j \in \mathbb{Z}}}, \max(|Y_{j}^n - Y_{j-1}^n|,|Y_{j+1}^n - Y_{j}^n|) > \epsilon \}
$
for a given $\epsilon > 0$. Following Karni \cite{karni2}, we will use $\epsilon = 0.05$ in practice. \\
\ \\
{\it{The quasi-conservative $p$-projection step with sampling ($t^{n+1-} \rightarrow t^{n+1}$)}} \\
The quasi-conservative $p$-projection step will be seen in the next section to properly compute 
large amplitude shock propagations. Localizing the averaging process of $p$ 
nevertheless prevents the method from keeping strictly uniform the velocity and pressure profiles of an isolated material 
front, see {\bf Test A} below. 
Indeed, note that since ${{\mathbb{Z}_p^{\varepsilon}}}$ is generally a {\it strict} 
subset of $\mathbb{Z}$ due to the numerical diffusion on $Y$ ({\it i.e.} ${{\mathbb{Z}_p^{\varepsilon}}} \subsetneq \mathbb{Z}$), 
an usual conservative treatment is still used on 
$\rho E$ as soon as $j$ not in ${{\mathbb{Z}_p^{\varepsilon}}}$.
This is sufficient to create pressure oscillations. 
In order to cure this problem, 
we propose to get rid of the numerical diffusion on $Y$ so as to 
enforce the non conservative treatment (\ref{newp}) 
across an isolated material front. 
This objective is achieved when replacing the conservative updating formula (\ref{moyL2phi1}) 
for $\rho_1$ with 
random sampling strategy applied to $Y$ (see also \cite{cchyp2006} and \cite{chalonsgoatin} for similar ideas). More precisely, 
we consider an equidistributed random sequence $(a_n)$ in $(0,1)$ 
(following Collela \cite{collela1}, we take in practice the celebrated van der Corput sequence), 
define $x^*_j = x_{j-1/2} + a_{n+1} \Delta t$ for all $j \in \mathbb{Z}$ and set 
$$
\left\{
\begin{array}{rcl}
Y^{n+1}_{j} = Y^{n+1-}_{j-1} & \mbox{if} & x^*_j \leq x_{j-1/2}^*, \\
Y^{n+1}_{j} = Y^{n+1-}_{j} & \mbox{if} & x_{j-1/2}^* \leq x^*_j \leq x_{j+1/2}^*, \\
Y^{n+1}_{j} = Y^{n+1-}_{j+1} & \mbox{if} & x^*_j \geq x_{j+1/2}^*.
\end{array}
\right.
$$
Then, we set $(\rho_1)^{n+1}_{j} = \rho^{n+1}_{j} Y^{n+1}_{j}$ so that the conservation of $\rho_1$ 
now holds only 
statistically.

\section{Numerical results} \label{sec:numexp}
We propose two numerical experiments with $\gamma_1=1.4$ and $\gamma_2=1.6$ associated with a Riemann initial data.
The left and right vectors ${\bf v}$ are denoted ${\bf v}_L$ and ${\bf v}_R$ and the initial discontinuity is 
at $x=0.5$.
In the first simulation ({\textbf{Test A}}), we consider the propagation of an isolated material interface 
with ${\bf v}_L=(1,1,1,1)$ and ${\bf v}_R=(0.1,0.1,1,0)$. We take $\Delta x = 0.005$ and plot the solutions 
at time $t=0.15$.
The second simulation ({\textbf{Test B}}) develops a strong shock due to a large initial pressure ratio. More precisely, 
we choose ${\bf v}_L=(1,0,500,1)$ and ${\bf v}_R=(1,0,0.2,0)$. We take $\Delta x = 0.00125$ and plot the solutions 
at time $t=0.008$. \\
We observe that the Two-Flux method and the Lagrangian methods are in agreement 
with the exact solutions and give similar results. As expected, note that the Lagrangian approach without sampling 
does not strictly maintain uniform the pressure profile for {\bf Test A}. Note also that the mass fraction $Y$ is sharp 
when a random sampling is used. At last, the relative conservation error on $\rho E$ (see for instance \cite{chalonsgoatin} for 
more details) for the Lagrangian approach with random sampling is actually less important and swings around 0.2\% only.

\begin{figure}[h!]
\begin{minipage}[b]{0.5\linewidth}
\centering \epsfig{file=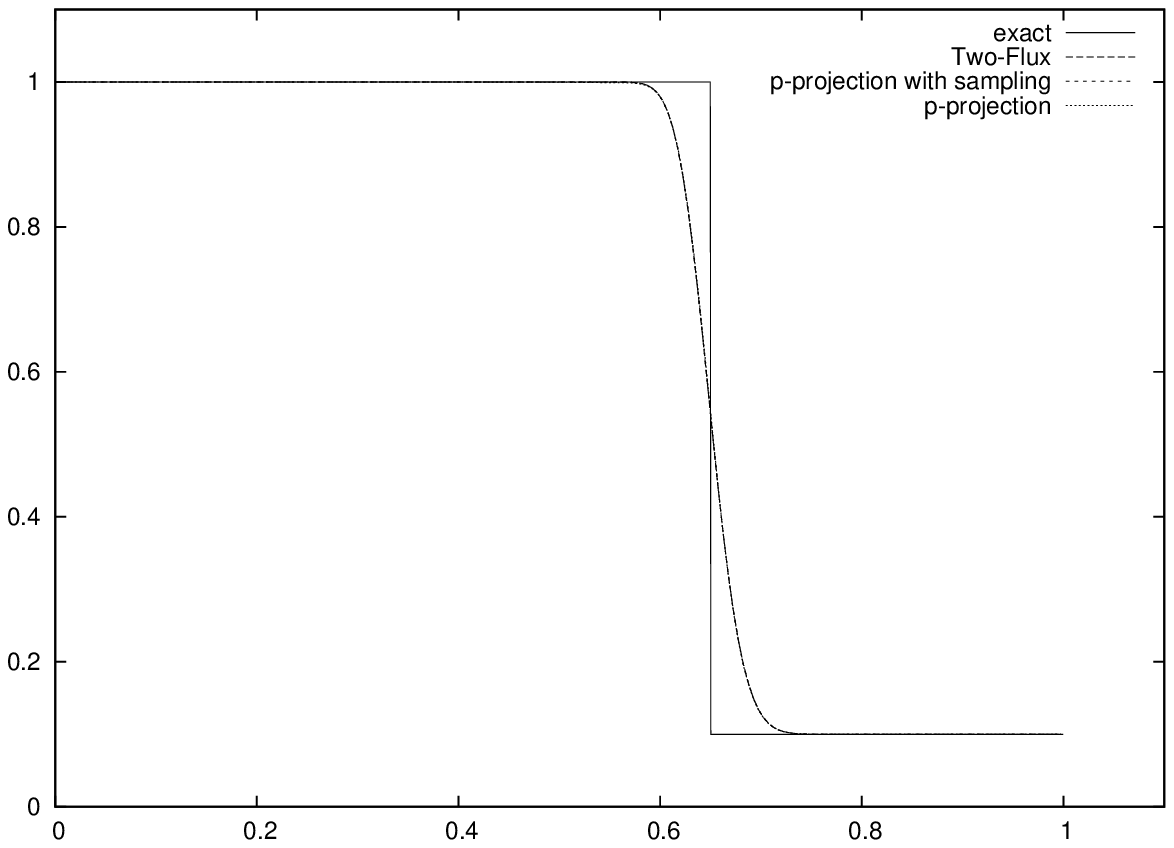,width=0.95\linewidth}
\vspace{0.3cm}
\caption{$\rho$ ({\bf Test A})}
\end{minipage}%
\begin{minipage}[b]{0.5\linewidth}
\centering \epsfig{file=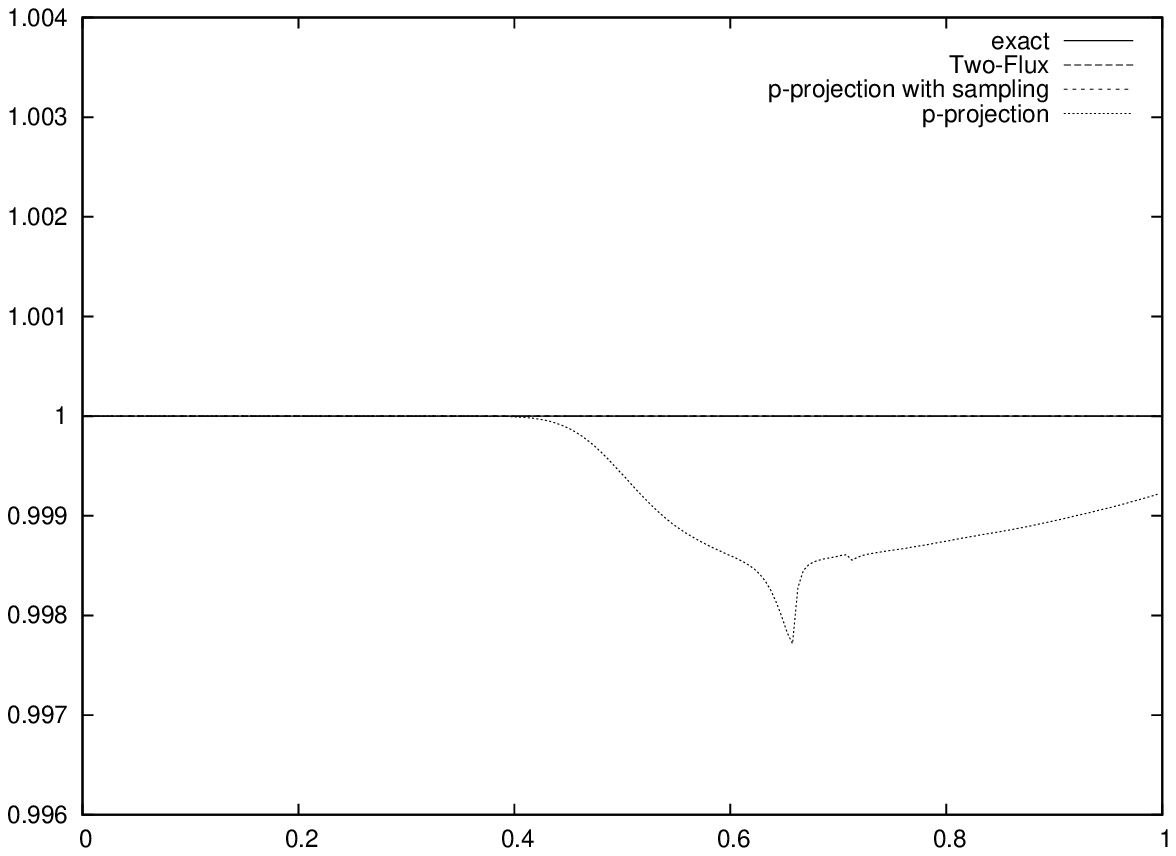,width=0.95\linewidth}
\vspace{0.3cm}
\caption{$p$ ({\bf Test A})}
\end{minipage}
\end{figure}
\begin{figure}[h!]
\begin{minipage}[b]{0.5\linewidth}
\centering \epsfig{file=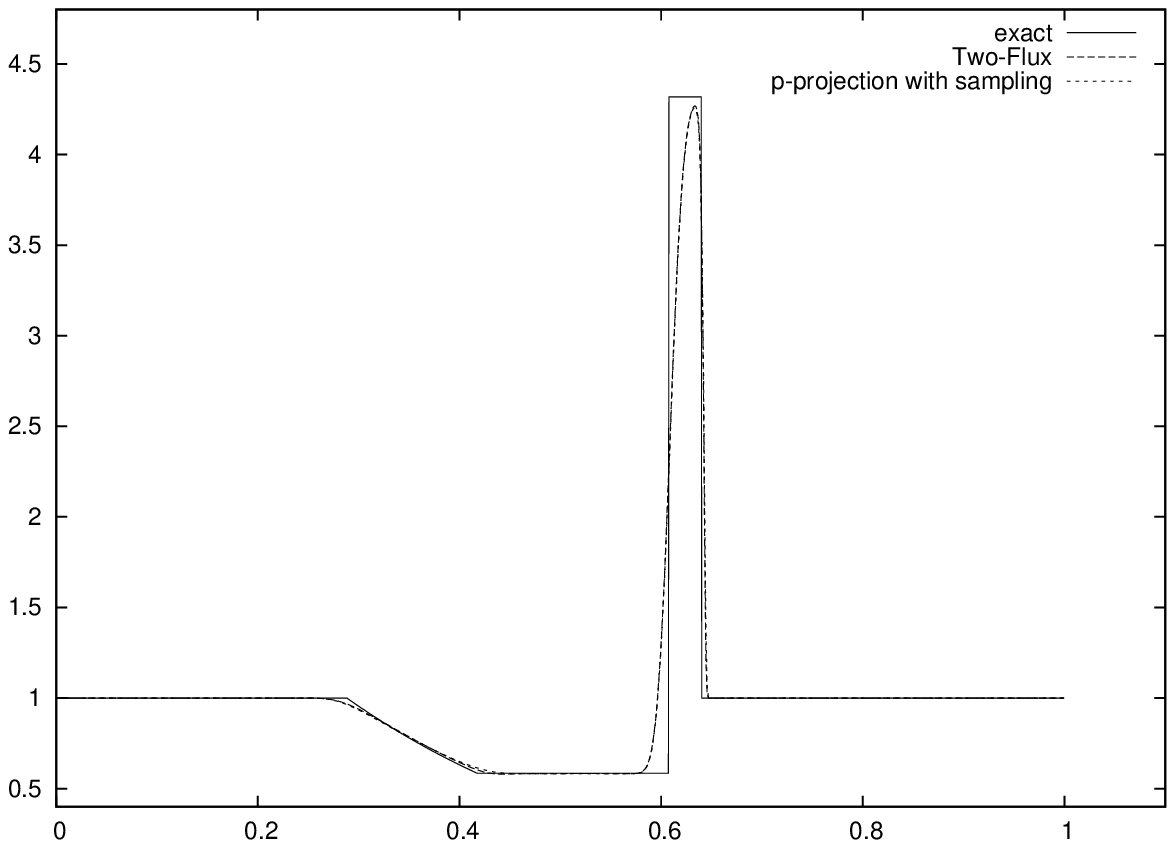,width=0.95\linewidth}
\vspace{0.3cm}
\caption{$\rho$ ({\bf Test B})}
\end{minipage}%
\begin{minipage}[b]{0.5\linewidth}
\centering \epsfig{file=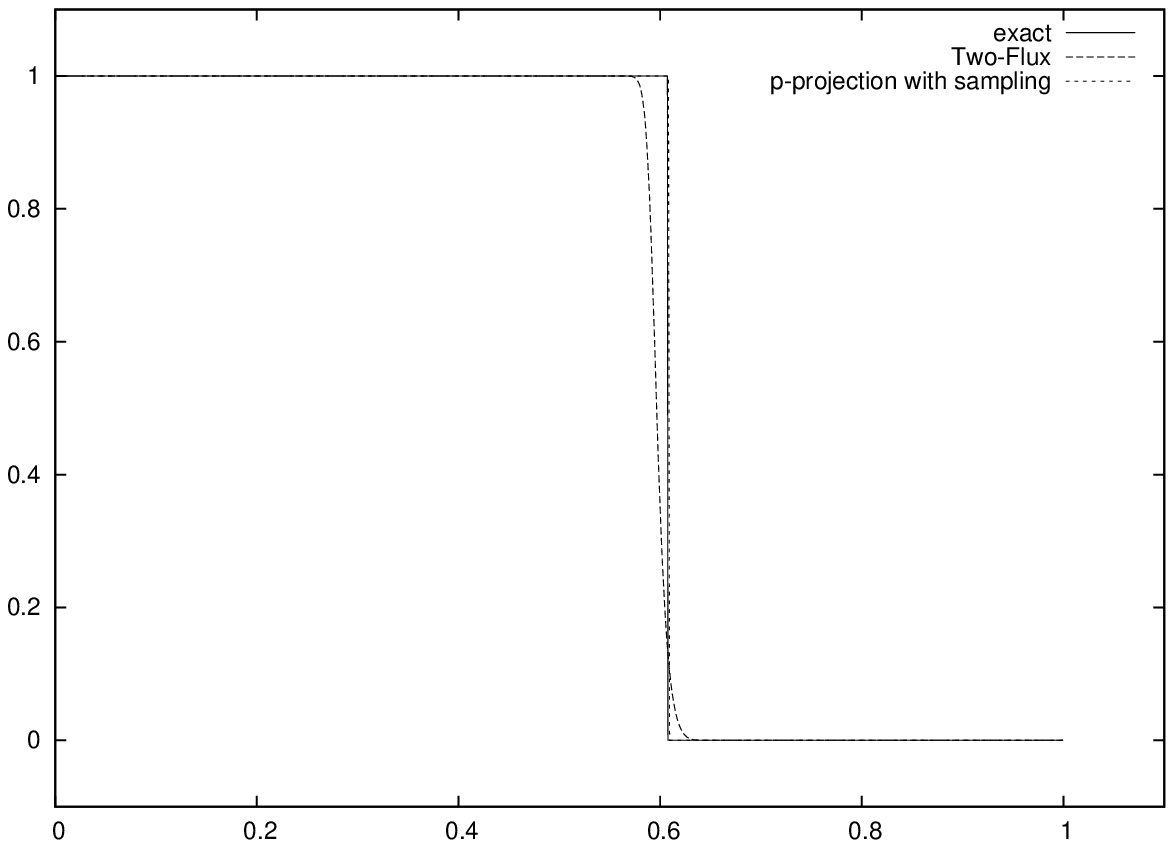,width=0.95\linewidth}
\vspace{0.3cm}
\caption{$Y$ ({\bf Test B})}
\end{minipage}
\end{figure}
\begin{figure}[h!]
\begin{minipage}[b]{0.5\linewidth}
\centering \epsfig{file=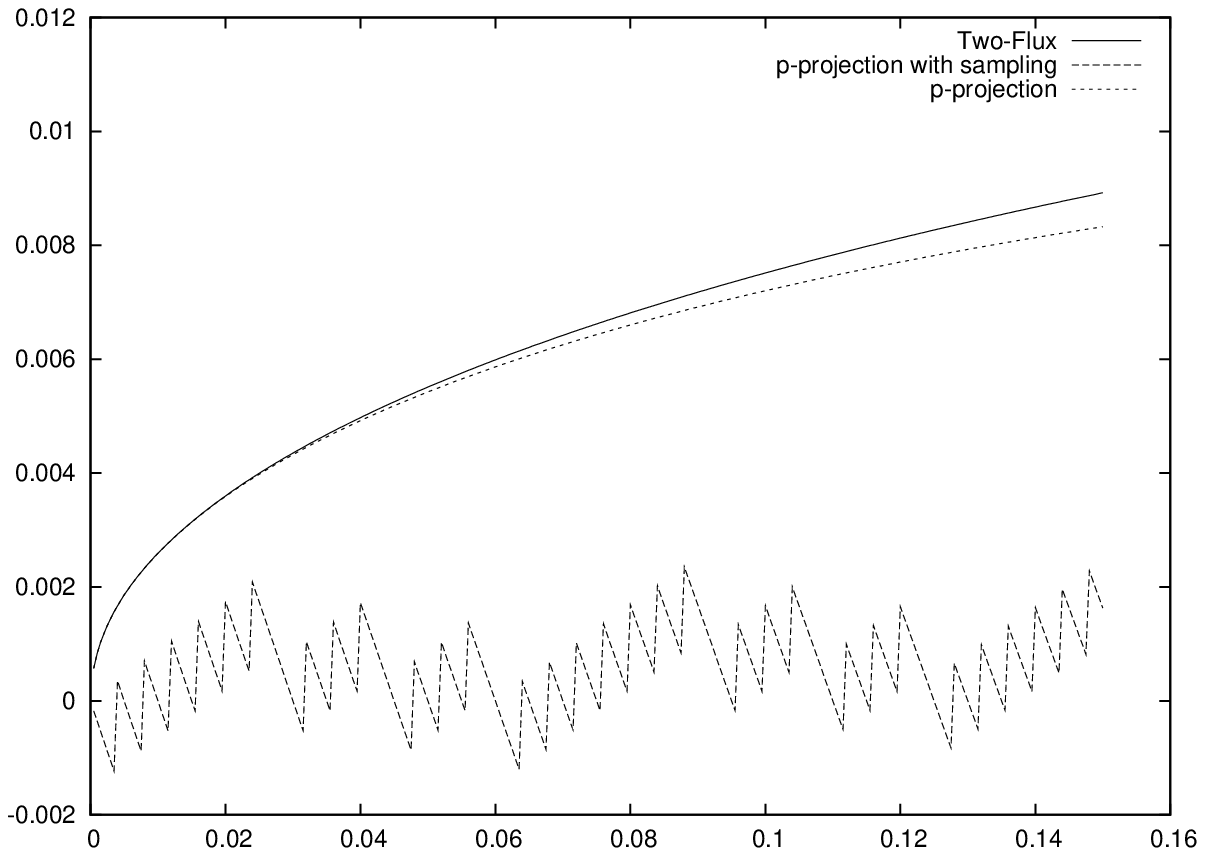,width=0.95\linewidth}
\vspace{0.3cm}
\caption{Conservation errors ({\bf A})}
\end{minipage}%
\begin{minipage}[b]{0.5\linewidth}
\centering \epsfig{file=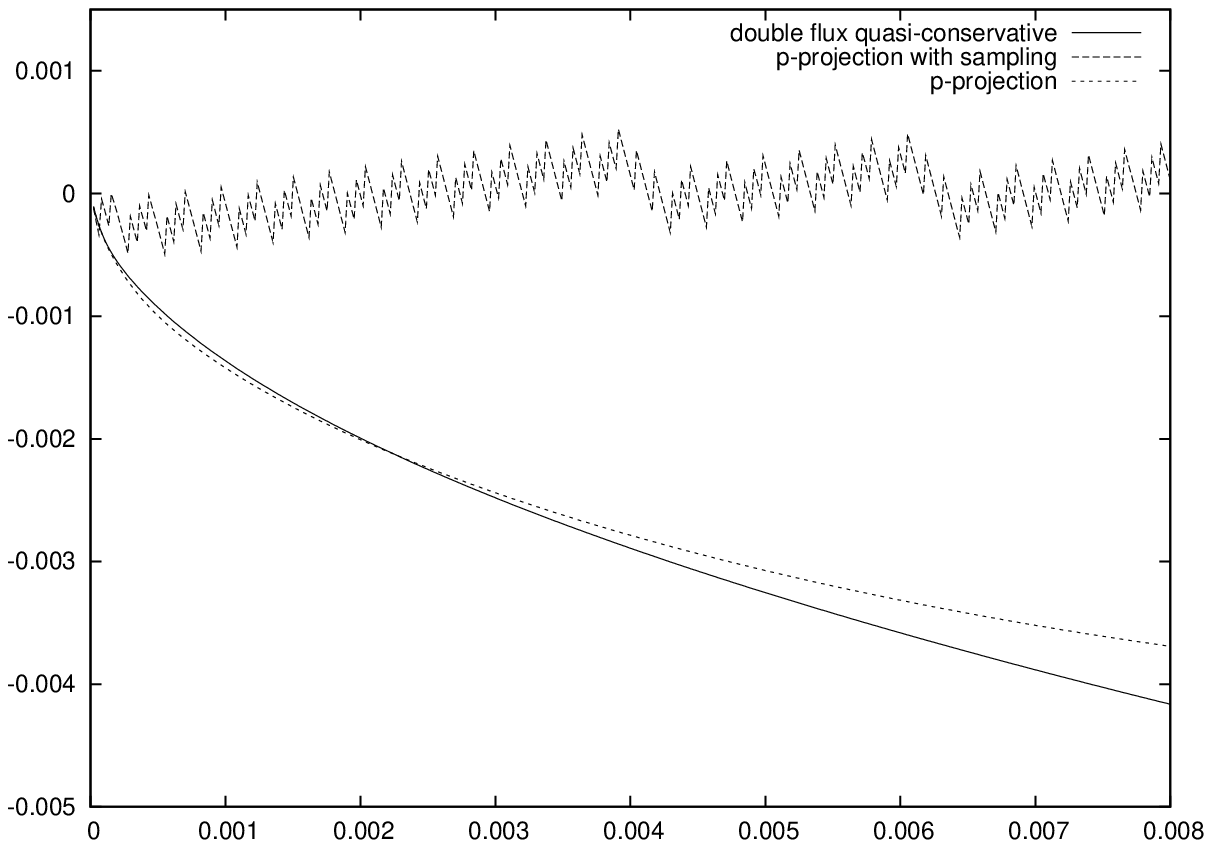,width=0.95\linewidth}
\vspace{0.3cm}
\caption{Conservation errors ({\bf B})}
\end{minipage}
\end{figure}

\section{Concluding remarks}
We have investigated a Lagrange-Projection approach 
for computing material fronts in multi-fluid models. We get similar 
results to the Two-Flux method \cite{ak1} with less important 
conservation errors on $\rho E$. Let us mention that other strategies, like for instance 
the one consisting in a local random sampling of ${\bf u}$ (instead of $Y$ 
only) in the Lagrangian step, have been investigated. The results 
are not reported here. \\
\ \\
\noindent {\bf Acknowledgements.} The authors are grateful for helpful discussions and exchanges with 
P. Helluy and F. Lagouti\`ere.

\end{document}